\magnification 1200
\def\R{{\rm I\kern-0.2em R\kern0.2em \kern-0.2em}}
\def\N{{\rm I\kern-0.2em N\kern0.2em \kern-0.2em}}
\def\P{{\rm I\kern-0.2em P\kern0.2em \kern-0.2em}}
\def\B{{\rm I\kern-0.2em B\kern0.2em \kern-0.2em}}
\def\Z{{\rm I\kern-0.2em Z\kern0.2em \kern-0.2em}}
\def\C{{\bf \rm C}\kern-.4em {\vrule height1.4ex width.08em depth-.04ex}\;}

\def\D{{\Delta}}

\def\z{{\zeta}}

\def\cH{{\cal H}}

\def\g{{\gamma}}
\font\ninerm=cmr8
\noindent  {\ninerm To appear in Math.\ Z.} 
\vskip 25mm
\centerline {\bf THE ARGUMENT PRINCIPLE AND HOLOMORPHIC EXTENDIBILITY}

\centerline {\bf TO FINITE RIEMANN SURFACES}
\vskip 4mm
\centerline{Josip Globevnik}
\vskip 4mm
{\noindent \ninerm ABSTRACT\ \  Let $M$ be a finite Riemann surface and let $A(M)$ be
the algebra of all continuous functions on $M\cup bM$ which are holomorphic on $M$. 
We prove that a 
continuous function $\Phi $ on $bM$ extends to a function in $A(M)$ if and only if 
 for every $f,\ g$ in $A(M)$ such that $f\Phi+g\not=0 $ on $bM$, the change 
 of argument of $f\Phi +g$ is nonnegative } 
\vskip 6mm

\bf 1.\ Introduction and the main result
\vskip 2mm
\rm
Let $M$ be a finite Riemann surface, that is, a region in a Riemann surface whose boundary
$bM$
consists of finitely many pairwise disjoint real analytic simple closed curves such that 
$M\cup bM$ is compact. We give $bM$ the standard orientation. Denote by $A(M)$ the algebra 
of all continuous functions on $M\cup bM$ which
are holomorphic on $M$. We denote by $\Delta $ the open unit disc in $\C$. 

The present paper deals with the problem of characterizing the boundary values of functions 
in $A(M)$ in terms 
of the argument principle. Our main result is
\vskip 2mm
\noindent \bf THEOREM 1.1\ \it Let $M$ be a finite Riemann surface. A continuous function 
$\Phi $ on $bM$ 
extends holomorphically through $M$ if and only if for every $f, g$ in $A(M)$ such that 
$f\Phi +g\not= 0$ on $bM$, the change of argument of $f\Phi + g$ along $bM$ is nonnegative. \rm
\vskip 2mm
In [Gl1] it was shown that $\Phi $ extends holomorphically through $M$ if and only if for 
every polynomial $P$ with coefficients in $A(M)$ such that $P(\Phi )\not=0$ on $bM$, the
change of argument of $P(\Phi )$ along $bM$ is nonnegative. This is an easy consequence 
of the fact that 
the algebra  $A(M)|bM$ is a maximal subalgebra of $C(bM)$ [R]. Theorem 1.1 sharpens this
by saying that in this characterization 
it suffices to take the polynomials of degree one.

The only if part of the theorem is an obvious consequence of the argument principle. In fact, if $\Phi$ 
admits the extension $\tilde \Phi \in A(M)$ then the change of the argument of $f\Phi +g$ along $bM$ 
equals $2\pi$ times the number of zeros of $f\tilde\Phi+g$ on $M$.
\vskip 4mm
\bf 2.\ Preliminaries
\vskip 2mm \rm
Let $g\geq 0$ be the genus of $M$ and let $m$ be the number of boundary components. Let 
$\nu = 2g+m-1$ and 
let simple closed curves $\g _1,\ \g _2,\cdots \g _\nu \subset M$ form a canonical basis for
$M$ which, together with $\g _1^\prime ,\ \g _2^\prime ,\cdots \g _\nu^\prime 
$ forms a symmetric canonical basis for the double $M^\ast $ of $M$ a 
closed Riemann surface of genus
$2\nu$ [K]. 

Let $U\subset M^\ast $ be a small open connected neighbourhood of $M\cup bM$ and let $W\subset 
M^\ast \setminus U$ be an open set which does not meet any of the curves 
$\g _1^\prime ,\ \g _2^\prime ,\cdots \g _\nu^\prime $. It is known that there is a 
meromorphic differential on $M^\ast $ with 
simple poles all contained in $W$ and with prescribed periods along the 
curves of the 
symmetric canonical basis of $M^\ast $ [K, p.19]. In particular, for each $j,\ 1\leq j\leq \nu$,
there 
is such a differential whose period along $\g _j$ equals $i$ and whose periods along all 
$\g _k,\ k\not= j$, are $0$. On $U$ the real part of such a differential is the differential
of a single valued harmonic function $u_j$ on $U$ whose conjugate differential $\ast d u_j$ 
has period one 
along $\g _j$ and period $0$ along 
each $\g _k,\ k\not= j$. (See [\v CF] for a different construction of such functions ). It 
follows that given a complex valued harmonic function  $u$ on $M$ there are unique 
constants
$c_1(u), c_2(u),\cdots c_\nu(u)$ such that $u-\sum_{j=1}^\nu c_j(u)u_j$ has a single 
valued conjugate, that is, 
$$
u(z)-\sum_{j=1}^\nu c_j(u)u_j(z) = F(z)+\overline{G(z)}\ \ (z\in M)
$$
 where $F$ and $G$ are 
(single valued) holomorphic functions on $M$. The constants $c_j(u)$ 
are the periods of the conjugate differential $\ast du$ and it is easy to see that they
depend continuously on $u$ in the sup norm. 
If the function $u$ extends smoothly to $M\cup bM$ then the same holds for the 
functions $F$ and $G$, so 
in this case $F$ and $G$ belong to $A(M)$ [B, p.91].
\vskip 4mm
\bf 3.\ Functions with single valued conjugates \rm
\vskip 2mm
If $\Phi $ is a continuous function on $bM$ there is a unique continuous extension 
$\cH (\Phi )$ of
$\Phi $ to $M\cup bM$ which is harmonic on $M$. We shall say that $\Phi $ \it has a 
single valued conjugate \rm if $\cH (\Phi )$ has a single valued 
conjugate on $M$. If this is the case then $\cH (\Phi ) = F + \overline G $
where $F$ and 
$G$ are holomorphic functions on $M$. 
In the special case when $\Phi $ is smooth the functions $F$ and $G$ belong to $A(M)$.
\vskip 2mm
\noindent\bf
Proposition 3.1\ \it Let $\Phi $ be a continuous function on $bM$. There is a 
nonconstant function $g$, holomorphic on $U$, such that the function 
$z\mapsto g(z)\Phi (z) \ (z\in bM)$ has a single valued conjugate. \rm
\vskip 2mm
\noindent \bf Proof.\ \rm Let $\omega $ be a nonconstant holomorphic function 
on $U$. For each $j,\ 1\leq j\leq \nu +1$, there are constants 
$c_{j 1},c_{j2},\cdots , c_{j\nu } $ such that the function 
$$
\z \mapsto \omega (\z )^j\Phi (\z) + \sum _{k=1}^\nu c_{j k}u_k(\z )\ \ \ (\z \in bM)
$$
has a single valued conjugate. Since $\nu +1$ rows $[c_{j 1},c_{j2},\cdots , c_{j\nu }]
\ \ (1\leq j\leq \nu+1) $ are èinearly dependent there are numbers 
$\lambda _1, \lambda_2,\cdots \lambda _{\nu +1}$, not all equal to zero, such that 
$\sum _{j=1}^{\nu+1}\lambda _j [c_{j 1},c_{j2},\cdots , c_{j\nu }] = 
[0, 0, \cdots 0]$ which means that the function $\z\mapsto\lambda _1\omega(\z )\Phi (\z )
+$ $\cdots \lambda_{\nu+1}$ $
\omega (\z ) ^{\nu +1}\Phi (\z )$\ $(\z\in bM)$ has a single valued conjugate. Thus, there is a polynomial 
$P$ with complex coefficients of degree at least one such that 
the function $\z\mapsto P(\omega (\z))\Phi (\z)$\  $(\z\in bM)$ has a single valued conjugate.
Assume for a moment that there is a constant $c$ such that
$$
P(\omega (\z )) + c \equiv 0\ \ (\z \in bM).
\eqno (3.1)
$$
Since $P$ has degree at least one it follows that $P+c$ has degree at least
one so there are 
$k\geq 1$ and complex numbers $\alpha \not=0,\ \omega _1, \omega _2,\cdots ,\omega _k$
such that 
$P(z)+c= \alpha (z - \omega _1)\cdots (z -\omega _k)$ so (3.1) implies that 
$$
(\omega (\z )-\omega _1)\cdots (\omega (\z )-\omega _k)           \equiv 0\ \ \ (\z \in bM).
\eqno (3.2)
$$
Since $\omega $ is holomorphic and nonconstant on $U$ it follows that each factor
in (3.2) 
has at most finiteèly many zeros on $bM$. This shows that (3.2) and hence (3.1) is 
impossible. This proves that 
$g=P(\omega )$ is a nonconstant function, holomorphic on $U$ and such
that the function $\z\mapsto g(\z )\Phi (\z )\ (\z \in bM)$ has a single valued conjugate. 
The proof is complete. 
\vskip 2mm
\noindent\bf Proposition 3.2\ \it Let $\Phi $ be a continuous function on 
$bM$ which does
not extend holomorphically through $M$. Given a nonconstant 
holomorphic function $g$ on $U$ there is an $a\in M$ such that 
$\cH ((g-g(a)\Phi )(a)\not= 0$.    
\vskip 2mm
\noindent \bf Proof.\ \rm Suppose that     $\cH ((g-g(a)\Phi )(a)\equiv 0 \ \ (a\in M)$.
It follows that 
$\cH (g\Phi )(a)= g(a)\cH (\Phi )(a)\ (a\in M)$ which, in particular, implies that 
$a\mapsto g(a)\cH (\phi )(a)$ is harmonic on $M$.  
Let $a\in M$ and let 
$\varphi\colon\Delta\rightarrow M$ be a parametric disc, $\varphi (0)=a$. Then 
$\z\mapsto g(\varphi (\z ))\cH (\Phi )(\varphi (\z ))$, a product of  $g\circ\varphi $ a 
nonconstant holomorphic function 
on $\Delta $ and  $\cH (\Phi )\circ\varphi $, a harmonic function on $\Delta$  
is harmonic on $\Delta $ which implies 
[GL2] that $\cH (\Phi )\circ\varphi $ is holomorphic on $\Delta $. Thus, $\cH (\Phi )$ is
holomorphic in a neighbourhood of $a$ and since $
\cH (\Phi )$ is harmonic on $M$ it follows that $\cH (\Phi )$ is holomorphic on $M$ which 
is impossible since $\Phi $ 
does not extend holomorphically through $M$. This completes the proof. 
\vskip 4mm
\bf 4.\ The proof in the case when $\Phi$ has a single valued conjugate
\vskip 2mm
\noindent Proposition 4.1\ \it Suppose that $\Phi $ is a continuous function on $bM$ 
which has a single valued 
conjugate and which does not extend holomorphically through $M$. There are functions 
$P, Q \in A(M)$ such that 
$P\Phi + Q\not= 0$ on $bM$ and such that the change of argument of $P\Phi + Q$ along
$bM$ is negative. 
\vskip 2mm
\noindent \bf Proof.\  \rm By Proposition 3.1 there is a function $g$, holomorphic
and nonconstant on $U$, such that the function
$\z\mapsto g(\z )\Phi (\z )\ (\z\in bM)$ has a single valued conjugate. Since $\Phi $ has 
a single valued conjugate it follows that for each $a\in M$
the function $\z\mapsto [g(\z )-g(a)]\Phi (\z )\ (\z\in bM)$  has a single valued conjugate. 
Since $\Phi $ does not extend holomorphically through $M$,
Proposition 3.2 implies that $\cH ((g-g(a))\Phi )(a)\not= 0$ for some $a\in M$. 
With no loss of generality assume that

$$
\cH ((g-g(a))\Phi )(a) = 5\eta >0.
\eqno (4.1)
$$
To make the proof easier to understand we first show how we 
complete the proof in the special case when $\Phi $ is smooth. In this case 
there are  
$F, G \in A(M),\ F(a)=G(a)=0$,
such that 
$$
\cH ((g-g(a))\Phi ) = F+\overline G+5\eta .
\eqno (4.2)
$$
It follows that
$$
(g(\z )-g(a))\Phi (\z ) - F(\z )-G(\z )\in 5\eta + i\R\ \ (\z\in bM)
$$
which implies that
$$
(g(\z )-g(a))\Phi (\z ) - F(\z )-G(\z ) \not= 0\ \ (\z\in bM)
$$
and that the change of argument of the function
$\z\mapsto (g(\z )-g(a))\Phi (\z ) - F(\z )-G(\z )$ along $bM$ is zero. 

Let $p_a$ be a holomorphic function on $U$ whose only zero on U 
is a single zero at $a$ [BS, p.566]. Note that $P= (g-g(a))/p_a \in A(M)$ and, since $F(a)=G(a)=0$ 
it follows also that $Q = -(F+G)/p_a \in A(M)$. Since  $(g-g(a))\Phi - F-G = p_a(P\Phi + Q)$ on $bM$, 
the argument principle implies that the change of argument of $P\Phi + Q$ along $bM$ equals $-2\pi $. This 
completes the proof in the case when $\Phi $ is smooth. 

We now proceed to the proof for general continuous $\Phi $. Given a continuous function $\Psi $ on $M$ there 
are constants $c_1, c_2,\cdots c_\nu $ such that 
$\z\mapsto (g(\z )-g(a))\Psi (\z )+\sum_{j=1}^\nu c_j u_j(\z )\ (\z\in bM)$ has a 
single valued conjugate. The constants 
$c_j,\ 1 \leq j\leq \nu $ depend continuously on $\Psi $ and they all vanish if $\Psi = \Phi$. 
So there is a $\delta >0$ such that 
$$
\bigl\vert \sum_{j=1}^\nu c_ju_j(z)\bigr\vert <\eta\ \ (z\in M\cup bM)
\eqno (4.3)
$$
and
$$
\vert (g(z)-g(a))(\Psi (z)-\Phi (z))\vert <\eta\ \ (z\in bM)
\eqno (4.4)
$$
provided that $|\Psi - \Phi|<\delta $ on $bM$. Let $\Psi $ be a smooth function 
on $bM$ such that $|\Psi - \Phi| <\delta $ 
on $bM$. The function $(g-g(a))\Psi + \sum_{j=1}^\nu c_ju_j$ is smooth on $bM$ and 
has a single valued conjugate, so there are
$F, G\in A(M),\ \ F(a)=G(a)=0$, and a constant $\gamma $ such that 
$$
\cH ((g-g(a))\Psi )+ \sum_{j=1}^\nu c_ju_j  = F+\overline G +\gamma \hbox{\ on\ } M\cup bM
$$
so 
$$
\gamma = \cH ((g-g(a)\Psi )(a)+\sum_{j=1}^\nu c_ju_j (a).
$$
By the maximum principle (4.4) implies that  $|\cH((g-g(a))(\Psi -\Phi)(a)|<\eta $ so (4.2) and (4.3) imply that 
$$
|\gamma -5\eta|<2\eta 
\eqno (4.5)
$$
On $bM$ we have $(g-g(a))\Phi -F-G =$ $(g-g(a))(\Phi -\Psi)+(g-g(a))\Psi - F-G =$ $ 
(g-g(a))(\Phi -\Psi )+\gamma -\sum_{j=1}^\nu c_ju_j +G-\overline G$ which, by (4.3),
(4.4) and (4.5) implies 
that 
$$
\bigl\vert 5\eta - [(g-g(a))(\Phi-\Psi)+\gamma - \sum_{j=1}^\nu c_ju_j ] \bigr\vert <4\eta \hbox{\ on
\ }bM
$$
which implies that
$$
(g(\z )-g(a))\Phi (\z )-F(\z )-G(\z )\in [\eta, 9\eta] + i\R\ \ (\z\in bM).
$$
So  $(g-g(a))\Phi-F-G          \not= 0$ on $bM$ and the change of argument of 
$(g-g(a))\Phi-F-G$ along $bM$ is zero. We now complete the proof in the same 
way as in the case when $\Phi $ was smooth. The proof is complete.  
\vskip 4mm
\bf 5.\ Completion of the proof of Theorem 1.1
\vskip 2mm
\noindent\bf
Proposition 5.1\ \it Let $f,g\in A(\Delta ) $. Assume that $g\not= 0$ on $b\D\setminus\{ 1\},\ 
g(1)=0$, and that $g$ extends holomorphically into a neighbourhood of $1$. Assume that the
function $\z\mapsto f(\z )/g(\z )\ (\z\in b\D \setminus\{ 1\} ) $
extends continuously to $b\D$ . Then it extends to a function from $A(\D )$. 
\vskip 2mm
\noindent\bf Proof.\ \rm We first prove the proposition in the case when $g(\z )=\z -1$. 
The following proof was shown to the author by Miran \v Cerne. Denote by $p$ the 
continuous extension of $\z\mapsto f(\z)/(\z -1)\  
(\z\in b\D \setminus \{ 1\})$ to $b\D $ and let $\sum a_ne^{in\theta}$ be the Fourier series 
of $p$. Let $\sum b_n e^{in\theta}$ be the Fourier series of $\z\mapsto 
f(\z ) = (\z -1)p(\z )$. Clearly $b_n= a_{n-1}-a_n\ (n\in Z)$. Since $f$ belongs to $A(\D )$ we have $b_n =0\ (n<0)$ 
which implies that $a_{n-1}=a_n\ (n<0)$. Since $\lim _{n\rightarrow\pm\infty}a_n = 0$ it follows that $a_n=0 \ (n<0)$
which implies that $p$ extends to a function in $A(\D )$. This completes the proof in
the special case when $g(\z )= \z -1$.  The repeated application of the special case 
proves the proposition in the case when $g(\z ) = (\z -1)^n,\ n\in\N$. In the general case write 
$g(\z )= (\z -1)^nh(\z )$ where $h\in A(\D ),\ h(\z )\not= 0 \ (\z \in \overline\D ) $ and $n\in \N$. The 
preceding discussion implies that the function $\z\mapsto (f/h)(\z )/(\z-1)^n\ 
(\z\in b\D\setminus\{ 1\})$ extends to a function from $A(\D )$ which shows that
 $\z\mapsto f(\z )/g(\z )\ (\z\in b\D\setminus\{ 1\})$ extends to a function from $A(\D )$. 
 This completes the proof. 
 \vskip 2mm
 \noindent\bf Proposition 5.2\ \it Let $\Phi$ be a continuous function on $bM$ 
 and let $g$ be a holomorphic function on 
 $U$, $g\not\equiv 0$. Assume that the function $\z\mapsto g(\z )\Phi (\z )$ 
 extends holomorphically 
 through $M$. Then there is a holomorphic function $h$ on $U$ without a
 zero on $bM$ such that 
 the function $\z\mapsto h(\z )\Phi (\z )$ extends holomorphically through $M$. 
 \vskip 2mm
 \noindent\bf Proof.\ \rm By the assumption there is a function $G\in A(M)$ such that 
 $g\Phi = G$ on $bM$. If $g$ is a constant there is nothing to prove 
 so assume that $g$ is not a constant and $g(b)=0$ for some $b\in bM$. The 
 point $b$ is an isolated zero of $g$, let its degree be $k$. Let $p_b$ be a 
 holomorphic function on $U$ whose only zero on $U$ is a single zero at $b$. The function 
 $g_1= g/p_b^k$ is holomorphic on $U$ and has no zero at $b$ so the function $g_1\Phi $
 is continuous on $bM$. This means that the function $\z\mapsto G(\z ) /p_b(\z )^k$ 
 $(\z\in bM\setminus \{b\})$ extends continuously to $bM$. Using Proposition 5.1 we see that            
the function $\z\mapsto G(\z ) /p_b(\z )^k$ extends to a function $G_1\in A(M)$. So, 
on $bM$ we have $g_1\Phi = G_1$ where $G_1 \in A(M)$ and where $g_1$ is holomorphic on 
$U$ and has the same zeros as $g$ except at $b$ where $g_1$ is different from $0$. 
Repeating the process at each zero of $g$ contained in $bM$ (there are only finitely 
many of these) we arrive at a function $h$ holomorphic on $U$ with no zero on $bM$ such that 
$\z\mapsto h(\z )\Phi (\z )$ extends to a function from $A(M)$. This completes the proof. 
\vskip 2mm
\noindent\bf Proof of Theorem 1.1 continued.\ \rm Suppose that $\Phi $ does not extend 
holomorphically through $M$. 
In the case when $\Phi $ has a single valued conjugate Proposition 4.1 implies that 
there are $P, Q\in A(M)$ such that $P\Phi+Q\not=0$ 
on $bM$ and such that the change of argument of $P\Phi + Q$ along $bM$ is negative. 
Suppose now that $\Phi $ does not 
have a single valued conjugate. By Proposition 3.1 there is a nonconstant function 
$g$ holomorphic on $U$ such that $(g|bM)\Phi $ has a single valued conjugate. If  
$(g|bM)\Phi $ 
does not extend holomorphically through $M$ Proposition 4.1 implies that there are $P_1, 
Q \in A(M)$ such that $P_1g\Phi + Q\not= 0$ on $bM$ and such that the change of argument of $
P_1g\Phi + Q$ along $bM$ is negative. Putting $P= P_1g$ completes the proof in the case when 
$(g|bM)\Phi $ does not extend holomorphically through $M$. Suppose now that there is a 
$G\in A(M)$ such that 
$g(\z )\Phi (\z ) = G(\z )\ \ (\z\in bM)$. By Proposition 5.2 there are $H\in A(M)$ and 
a holomorphic function $h$ on $U$ having no zero on $bM$ such that $h(\z )\Phi (\z )= H(\z )\ 
(\z\in bM)$. Dividing both sides with powers of functions $p_b,\ b\in M$ where 
$p_b$ is a holomorphic function on $U$ whose only zero on $U$ is a single zero at $b$,
we may 
assume with no loss of generality that $h$ and $H$ have no common zero on $M$. Since $h$
has no zero on $bM$ 
and since $\Phi $ does not extend holomorphically through $M$ it follows that there is 
an $a\in M$ which 
is a zero of $h$ and not a zero of $H$. On $bM$ we have
$$
{h\over {p_a}}\Phi = {H\over{p_a}} ={{H-H(a)}\over{p_a}} + {{H(a)}\over{p_a}}
$$ 
where $H(a)\not= 0$. The functions $P=h/p_a$ and $Q= -(H-H(a))/p_a$ both belong to $A(M)$ 
and we have $P\Phi + Q = H(a)/p_a $ on $bM$ which implies that $P\Phi+Q\not =0$ on $bM$ and, 
by the argument principle, the change of argument of $P\Phi + Q$ along $bM$ is negative. 
The proof of Theorem 1.1 is complete. 
\vskip 4mm
\noindent\bf 6.\ A remark \rm
\vskip 2mm
If $M$ is a finitely connected domain in $\C$ then 
$\Phi\in C(bM)$ extends holomorphically through $M$ provided that for each $Q\in A(M)$ 
such that 
$\Phi+Q\not=0 $ on $bM$, the change of argument of $\Phi + Q$ along $bM$ is nonnegative 
(that is, 
if $M$ is a finitely connected domain in the plane then in Theorem 1.1 it suffices to take 
$P\equiv 1$). The question whether the same holds for general finite Riemann surfaces remains open. 

\vskip 6mm
\noindent \bf Acknowledgements \rm\ The author is indebted to Miran \v Cerne for the 
proof of Proposition 5.1. 
A major part of the work whose results are presented here was done in the Fall of 2004 during the 
author's visit 
at the Department of Mathematics, University of California, San Diego.
The author
is indebted 
to Salah Baouendi, Peter Ebenfelt and Linda Rothschild for making the very pleasant 
visit possible.  

This work was supported 
in part by the Ministry of Higher Education, Science and Technology of Slovenia  
through the research program Analysis and Geometry, Contract No.\ P1-0291. 
\vfill
\eject

\centerline{\bf REFERENCES} 
\vskip 5mm
\noindent [AW]\ H.\ Alexander and J.\ Wermer: Linking 
numbers and boundaries of varieties. 

\noindent Ann.\ Math.\ 151 (2000) 125-150
\vskip 2mm
\noindent [BS]\ \ H.\ Behnke, F.\ Sommer: \it Theorie der analitischen Funktionen 
einer komplexen
Ver\" an\-derlichen. \rm

\noindent Springer, Berlin-Gottingen-Hiedelberg 1955
\vskip 2mm
\noindent [B]\ \ S.\ Bell: \it The Cauchy transform, Potential Theory, and 
Conformal Mapping. \rm CRC Press, 
Boca Raton, 1992
\vskip 2mm
\noindent [\v CF]\ \ M.\ \v Cerne, M.\ Flores:\ Generalized Ahlfors functions. 

\noindent To appear in Trans.\ Amer.\ Math.\ Soc.
\vskip 2mm
\noindent [Gl1]\ \ J.\ Globevnik:\ Holomorphic extendibility and the argument principle.

\noindent To appear in "Complex Analysis and Dynamical Systems II (Proceedings
  of a conference held in honor of Professor Lawrence Zalcman's sixtieth
  birthday in Nahariya, Israel, June 9-12, 2003)", Contemp.\ Math. 
  [http://arxiv.org/abs/math.CV/0403446]
\vskip 2mm
\noindent [Gl2]\ \ J.\ Globevnik:\ The argument principle and holomorphic extendibility.

\noindent Journ.\ d'Analyse.\ Math.\ 94 (2004) 385-395
\vskip 2mm
\noindent [K]\ \   W.\ Koppelman:\ The Riemann-Hilbert problem for finite Riemann surfaces. 

\noindent Comm.\ Pure Appl.\ Math.\ 12 (1959) 13-35
\vskip 2mm
\noindent [R]\ \ H.\ L.\ Royden:\ The boundary values of analytic and harmonic functions.

\noindent Math.\ Z.\ 78 (1962) 1-24
\vskip 2mm
\noindent [S]\ \  G.\ Springer: \ \it Introduction to Riemann Surfaces.\rm 

\noindent Addison-Wesley, 1957
\vskip 2mm
\noindent [W] J.\ Wermer: The argument principle and boundaries of analytic varieties.

\noindent Oper. Theory Adv. Appl., 127, Birkhauser, Basel, 2001, 639-659
\vskip 10mm

\noindent 
Institute of Mathematics, Physics and Mechanics

\noindent 
University of Ljubljana

\noindent Ljubljana, Slovenia

\noindent josip.globevnik@fmf.uni-lj.si

\end